\documentclass[12pt]{article}

\usepackage{amsmath, epsfig, cite}
\usepackage{amssymb}
\usepackage{amsfonts}
\usepackage{latexsym}
\usepackage{graphicx}
\usepackage{color}
\usepackage{ifpdf}
\usepackage{CJK}

\newtheorem{thm}{Theorem}[section]

\newtheorem{cor}[thm]{Corollary}

\newtheorem{lem}[thm]{Lemma}
\newtheorem{conj}[thm]{Conjecture}

\numberwithin{equation}{section}

\makeatletter \@addtoreset{equation}{section} \makeatother

\begin{document}
\rule{0cm}{1cm}

\begin{center}
{\Large\bf On a relation between the Szeged index\\[2mm] and the Wiener
index for bipartite graphs }
\end{center}

\begin{center}
{\small Lily Chen, Xueliang Li, Mengmeng Liu\\
Center for Combinatorics, LPMC-TJKLC\\
Nankai University, Tianjin 300071, China\\
Email:  lily60612@126.com, lxl@nankai.edu.cn, liumm05@163.com}
\end{center}

\begin{center}
\begin{minipage}{120mm}
\begin{center}
{\bf Abstract}
\end{center}

{\small The Wiener index $W(G)$ of a graph $G$ is the sum of the
distances between all pairs of vertices in the graph. The Szeged
index $Sz(G)$ of a graph $G$ is defined as $Sz(G)=\sum_{e=uv \in
E}n_u(e)n_v(e)$ where $n_u(e)$ and $n_v(e)$ are, respectively, the
number of vertices of $G$ lying closer to vertex $u$ than to vertex
$v$ and the number of vertices of $G$ lying closer to vertex $v$
than to vertex $u$. Hansen used the computer programm AutoGraphiX
and made the following conjecture about the Szeged index and the
Wiener index for a bipartite connected graph $G$ with $n \geq 4$
vertices and $m \geq n$ edges:
$$
Sz(G)-W(G) \geq 4n-8.
$$
Moreover the bound is best possible as shown by the graph composed
of a cycle on $4$ vertices $C_4$ and a tree $T$ on $n-3$ vertices
sharing a single vertex. This paper is to give a confirmative proof
to this conjecture. }

\vskip 3mm


\vskip 3mm


\end{minipage}
\end{center}

\section{Introduction}

All graphs considered in this paper are finite, undirected and
simple. We refer the readers to \cite{bm} for terminology and
notation. Let $G$ be a connected graph with vertex set $V$ and edge
set $E$. For $u,v \in V$, $d(u,v)$ denotes the {\it distance}
between $u$ and $v$. The {\it Wiener index} of $G$ is defined as
$$
W(G)=\displaystyle\sum_{\{u,v\}\subseteq V} d(u,v).
$$
This topological index has been extensively studied in the
mathematical literature; see, e.g., \cite{GSM,GYLL}. Let $e=uv$ be
an edge of $G$, and define three sets as follows:
$$
N_u(e) = \{w \in V: d(u,w)< d(v,w)\},
$$
$$
N_v(e) = \{w \in V: d(v,w)< d(u,w)\},
$$
$$
N_0(e) = \{w \in V: d(u,w)=d(v,w)\}.
$$
Thus, $\{N_u(e),N_v(e),N_0(e)\}$ is a partition of the vertices of
$G$ respect to $e$. The number of vertices of $N_u(e)$, $N_v(e)$ and
$N_0(e)$ are denoted by $n_u(e)$, $n_v(e)$ and $n_0(e)$,
respectively. A long time known property of the Wiener index is the
formula \cite{GP,W}:
$$
W(G) = \displaystyle\sum_{e=uv \in E} n_u(e) n_v(e),
$$
which is applicable for trees. Using the above formula, Gutman
\cite{G} introduced a graph invariant, named as the {\it Szeged
index} as an extension of the Wiener index and defined by
$$
Sz(G) = \displaystyle\sum_{e=uv \in E}n_u(e) n_v(e).
$$
Randi\'c \cite{R} observed that the Szeged index does not take into
account the contributions of the vertices at equal distances from
the endpoints of an edge, and so he conceived a modified version of
the Szeged index which is named as the {\it revised Szeged index}.
The revised Szeged index of a connected graph $G$ is defined as
$$
Sz^*(G) = \displaystyle\sum_{e=uv \in E}\left(n_u(e)+
\frac{n_0(e)}{2}\right)\left(n_v(e)+ \frac{n_0(e)}{2}\right).
$$

Some properties and applications of the Szeged index and the revised
Szeged index have been reported in \cite{AH,LL,PR,PZ,XZ}.

In \cite{auto}, Hansen used the computer programm AutoGraphiX and
made the following conjectures:

\begin{conj}\label{conj1}
Let $G$ be a bipartite connected graph with $n \geq 4$ vertices and
$m \geq n$ edges. Then
$$
Sz(G)-W(G) \geq 4n-8.
$$
Moreover the bound is best possible as shown by the graph composed
of a cycle on $4$ vertices $C_4$ and a tree $T$ on $n-3$ vertices
sharing a single vertex.
\end{conj}

\begin{conj}\label{conj}
Let $G$ be a bipartite connected graph with $n \geq 4$ vertices and
$m \geq n$ edges. Then
$$
Sz^*(G)-W(G) \geq 4n-8.
$$
Moreover the bound is best possible as shown by the graph composed
of a cycle on $4$ vertices $C_4$ and a tree $T$ on $n-3$ vertices
sharing a single vertex.
\end{conj}

It is easy to see that $Sz^*(G)=Sz(G)=W(G)$ if $G$ is a tree, which
means $m=n-1$. So, the second conjecture considers graphs with $m
\geq n$.

This paper is to give confirmative proofs to the two conjectures. In
fact, if $G$ is a bipartite graph, then $Sz^*(G)=Sz(G)$. Therefore,
if we give a proof to Conjecture \ref{conj1}, then Conjecture
\ref{conj} follows immediately.

\section{Main results}

In \cite{SGB}, Gutman gave another expression for the Szeged index:
$$
Sz(G)= \displaystyle\sum_{e=uv \in E}n_u(e) n_v(e)=
\displaystyle\sum_{e=uv \in E}\displaystyle\sum_{\{x,y\} \subseteq
V}\mu_{x,y}(e)
$$
where $\mu_{x,y}(e)$, interpreted as contribution of the vertex pair
$x$ and $y$ to the product $n_u(e)n_v(e)$, is defined as follows:

$$
\mu_{x,y}(e) =  \left\{
\begin{array}{ll}
1,& \mbox {if}
\begin{cases} d(x,u)<d(x,v) \text{ and } d(y,v)<d(y,u), \\ \text{or} \\
d(x,v)<d(x,u)\text{ and }\ d(y,u)<d(y,v),\\
\end{cases}\\
0, & \mbox{otherwise}.
\end{array}
\right.
$$

We first show that for a $2$-connected bipartite graph Conjecture
\ref{conj1} is true.

\begin{lem}\label{lem1}
Let $G$ be a 2-connected bipartite graph of order $n \geq 4$. Then
$$Sz(G)-W(G) \geq 4n-8$$
with equality if and only if $G=C_4$.
\end{lem}

\begin{pf}
From above expressions, we know that
\begin{eqnarray*}
Sz(G)-W(G) & = & \displaystyle \sum_{\{x,y\} \subseteq V}
\displaystyle \sum_{e \in E}\mu_{x,y}(e)-
\displaystyle \sum_{\{x,y\} \subseteq V}d(x,y)\\
&=& \displaystyle \sum_{\{x,y\} \subseteq V}\left( \displaystyle
\sum_{e \in E}\mu_{x,y}(e)-d(x,y)\right).
\end{eqnarray*}

\noindent {\bf Claim:} For every pair $x,y \in V$, we have
$$
\displaystyle \sum_{e \in E}\mu_{x,y}(e)-d(x,y) \geq 1.
$$

In fact, if $xy \in E$, that is $d(x,y)=1$, then we can find a
shortest cycle $C$ containing $x$ and $y$ since $G$ is
$2$-connected. Then, $G[C]$ has no chord. Since $G$ is bipartite,
the length of $C$ is even. There is an edge $e'$ which is the
antipodal edge of $e=xy$ in $C$. It is easy to check that
$\mu_{x,y}(e')=\mu_{x,y}(e)=1$. So the claim is true.

If $d(x,y) \geq 2$, let $P_1$ be a shortest path from $x$ to $y$ and
$P_2$ be a second shortest path from $x$ to $y$, that is, $P_2 \neq
P_1$ and $|P_2|=$ min $\{|P|| P $ is a path from $x$ to $y$ and
$P\neq P_1\}$. Since $G$ is $2$-connected, $P_2$ always exists. If
there are more than one path satisfying the condition, we choose
$P_2$ as a one having most common vertices with $P_1$.

If $E(P_1)\bigcap E(P_2)=\emptyset$, let $P_1 \bigcup P_2=C$, and
then $|E(P_2)|\geq |E(P_1)|$ and all the antipodal edges of $P_1$ in
$C$ makes $\mu_{x,y}(e)=1$. We also know that $\mu_{x,y}(e)=1$ for
all $e \in E(P_1)$. Hence, $\sum_{e \in E}\mu_{x,y}(e)-d(x,y) \geq
d(x,y)> 1.$

If $E(P_1)\bigcap E(P_2) \neq \emptyset$, then $P_1\triangle P_2=C$,
where $C$ is a cycle. Let $P'_i=P_i\bigcap C=x'P_iy'$. It is easy to
see that $|E(P'_2)| \geq |E(P'_1)|$, and the shortest path from $x$
(or $y$) to the vertex $v$ in $P'_2$ is $xP_2x'(yP_2y')$ together
with the shortest path from $x'(y')$ to $v$ in $C$; otherwise,
contrary to the choice of $P_2$. So, all the antipodal edges of
$P'_1$ in $C$ makes $\mu_{x,y}(e)=1$. We also know that
$\mu_{x,y}(e)=1$ for all $e\in E(P_1)$. Hence, $\sum_{e \in
E}\mu_{x,y}(e)=|E(P_1)|+ d(x',y')\geq d(x,y)+1,$ which proves the
claim.

Now let $C=v_1v_2\cdots v_pv_1$ be a shortest cycle in $G$, where
$p$ is even and $p \geq 4$. Actually, for every $e \in E(C)$ we have
that $\mu_{v_i,v_{\frac{p}{2}+i}}(e)=1$ for
$i=1,2,\cdots,\frac{p}{2}.$ Then $\sum_{e \in
E}\mu_{v_i,v_{\frac{p}{2}+i}}(e)=|C|=p,$ that is, $\sum_{e \in
E}\mu_{v_i,v_{\frac{p}{2}+i}}(e)-d(v_i,v_{\frac{p}{2}+i})=
\frac{p}{2} \geq 2.$ Combining with the claim, we have that
$$
Sz(G)-W(G) \geq {n \choose 2} + \frac{p}{2}
\left(\frac{p}{2}-1\right) \geq {n \choose 2} +2 \geq 4n-8.
$$
The last two equalities hold if and only if $p=4$, $n=4 \ \mbox{or}\
5$. If $n=4, p=4$, then $G$ is a $C_4$. If $n=5, p=4$, then $G$ is a
$K_{2,3}$, and in this case we can easily calculate that
$Sz(G)-W(G)>12$. Thus, the equality holds if and only if $G=C_4$.

\end{pf}
\begin{qed}
\end{qed}

Next we will complete the proof of Conjecture \ref{conj1} in
general.

\begin{thm}\label{thm2}
Let $G$ be a bipartite connected graph with $n \geq 4$ vertices and
$m \geq n$ edges. Then
$$
Sz(G)-W(G) \geq 4n-8.
$$
Moreover the bound is best possible as shown by the graph composed
of a cycle on $4$ vertices $C_4$ and a tree $T$ on $n-3$ vertices
sharing a single vertex.
\end{thm}

\begin{pf} We have proved that the conclusion is true for a $2$-connected
bipartite graph. Now suppose that $G$ is a connected bipartite graph
with blocks $B_1, B_2,\cdots, B_k$, where $k\geq 2$. Let
$|B_i|=n_i$. Then, $n_1+n_2+\cdots+n_k=n+k-1$. Since $m\geq n$ and
$G$ is bipartite, there exists at least one block, say $B_1$, such
that $n_1\geq 4$. Consider a pair $\{x,y\}\subseteq V$. We have the
following observations:

\noindent {\bf Obs.1:} $x, y\in B_i$, and $n_i\geq 4$. For every $e
\in B_j, j \neq i, \mu_{x,y}(e)=0$, combining with Lemma \ref{lem1},
we have that
$$\sum_{\{x,y\}\subseteq B_i}\left(\sum_{e \in E}\mu_{x,y}(e)-d(x,y)\right )
= \sum_{\{x,y\}\subseteq B_i}\left(\sum_{e \in
E(B_i)}\mu_{x,y}(e)-d(x,y)\right ) \geq 4n_i-8.$$

\noindent {\bf Obs.2:} $x, y\in B_i$, and $n_i= 2$. In this case,
$$\displaystyle\sum_{\{x,y\}\subseteq B_i}\left(\displaystyle \sum_{e \in
E}\mu_{x,y}(e)-d(x,y)\right ) =0=4n_i-8.$$

\noindent {\bf Obs.3:} $x\in B_1, y\in B_i, i\neq 1$. Let $P$ be a
shortest path from $x$ to $y$, and let $w_1,w_i$ be the cut vertices
in $B_1$ and $B_i$ such that every path from a vertex in $B_1$ to
$B_i$ must go through $w_1,w_i$. By the proof of Lemma \ref{lem1},
we can find an edge $e' \in E(B_1)\backslash E(P)$ such that
$\mu_{x,w_1}(e')=1$. Because every path from a vertex in $B_1$ to
$y$ must go through $w_1$, we have $\mu_{x,y}(e')=1$. We also know
that $\mu_{x,y}(e)=1$ for all $e \in E(P)$. Hence, $\sum_{e \in
E}\mu_{x,y}(e)-d(x,y) \geq 1.$

We are now in a position to show that for all $y \in B_i \backslash
\{w_i\}$, we can find a vertex $z \in B_1\backslash \{w_1\}$ such
that $ \sum_{e \in E}\mu_{z,y}(e)-d(z,y) \geq 2$. Since $B_1$ is
2-connected with $n_1\geq 4$, there is a cycle containing $w_1$. Let
$C$ be a shortest cycle containing $w_1$, say $C=v_1v_2\cdots
v_pv_1$, where $v_1=w_1,$ $p$ is even. Set $z=v_{\frac{p}{2}+1}$. By
the proof of Lemma \ref{lem1}, we have that $ \sum_{e \in
E(B_1)}\mu_{z,w_1}(e)-d(z,w_1) \geq \frac{p}{2}\geq 2$. It follows
that there are two edges $e',e''$ which are not in the shortest path
from $z$ to $w_1$ such that $\mu_{z,w_1}(e')=1, \mu_{z,w_1}(e'')=1$.
Thus, $\mu_{z,y}(e')=1, \mu_{z,y}(e'')=1$. Hence, $\sum_{e \in
E}\mu_{z,y}(e)-d(z,y) \geq 2$.

If we fix $B_i$, we obtain that
$$\displaystyle\sum_{\substack {x\in B_1\backslash\{w_1\}\\y\in
B_i\backslash\{w_i\}}} \left(\sum_{e \in
E}\mu_{x,y}(e)-d(x,y)\right) \geq
(n_1-1)(n_i-1)+(n_i-1)=n_1(n_i-1).$$

\noindent {\bf Obs.4:} $x\in B_i, y\in B_j, i\geq 2, j\geq 2, i\neq
j$. Let $P$ be a shortest path between $x$ and $y$. If $P$ passes
through a block $B_l$ with $n_l\geq 4$, and $|B_l\bigcap P|\geq2$,
then we have that $\displaystyle\sum_{e \in E}\mu_{x,y}(e)-d(x,y)
\geq 1$. Otherwise, $\displaystyle\sum_{e \in E}\mu_{x,y}(e)-d(x,y)
\geq 0.$ So,
$$\sum_{x\in B_i\backslash\{w_i\}, y\in B_j\backslash\{w_j\}}\left (
\sum_{e\in E}\mu_{x,y}(e)-d(x,y))\right)\geq 0.
$$
Equality holds if and only if $P$ passes through a block $B_l$ with
$n_l=2$ or $n_l \geq 4,$ and $|B_l \bigcap P|=1.$

From the above observations, we have that
\begin{eqnarray*}
& & Sz(G)-W(G) \\
& = & \displaystyle \sum_{\{x,y\} \subseteq V} \displaystyle \sum_{e
\in E}\mu_{x,y}(e)-
\displaystyle \sum_{\{x,y\} \subseteq V}d(x,y)\\
&=& \displaystyle \sum_{\{x,y\} \subseteq V}\left( \displaystyle
\sum_{e
\in E}\mu_{x,y}(e)-d(x,y)\right)\\
& =& \displaystyle\sum_{i=1}^{k}\sum_{\{x,y\} \subseteq B_i}
\left(\sum_{e \in E}\mu_{x,y}(e)-d(x,y)\right)
+\sum_{j=2}^{k}\sum_{\substack{x\in B_1\backslash\{w_1\}\\ y\in
B_j\backslash\{w_j\}}} \left(\sum_{e
\in E}\mu_{x,y}(e)-d(x,y)\right)\\
& & +\frac{1}{2}\sum_{\substack{i\neq j\\ i\neq 1,j\neq
1}}\sum_{\substack {x\in B_i\backslash\{w_i\}\\ y\in
B_j\backslash\{w_j\}}}\left (
\sum_{e\in E}\mu_{x,y}(e)-d(x,y))\right)\\
& \geq & \sum_{i=1}^{k}(4n_i-8)+n_1\sum_{j=2}^{k}(n_j-1)\\
& =& 4(n+k-1)-8k+n_1(n-n_1)\\
& =& 4n-4k-4+n_1(n-n_1).
\end{eqnarray*}

Since $n_1+n_2+\cdots+n_k=n+k-1$, $n_1\geq 4, n_i\geq 2,$ for $2\leq
i\leq k$, we have that $4\leq n_1\leq n-k+1,$ and $2\leq k\leq n-3.$

If $k\geq 5$, then $n_1(n-n_1)\geq 4(n-4)$. Thus,
$$4n-4k-4+n_1(n-n_1)\geq 8n-4k-20 \geq 8n-4(n-3)-20=4n-8.$$
Equality holds if and only if $n_1=4, n_2=n_3=\cdots=n_{n-3}=2$, and
$B_2, B_3, \cdots, B_{n-3}$ form a tree $T$ on $n-3$ vertices, which
shares a single vertex with $B_1$.

If $2\leq k\leq 4$, then $n_1(n-n_1)\geq (n-k+1)(k-1)$.

If $k=2$, then $4n-4k-4+(n-k+1)(k-1)=5n-13 \geq 4n-8.$ Equality
holds if and only if $n=5$, $G$ is a graph composed of a cycle on 4
vertices and a pendant edge.

If $k=3$, then $4n-4k-4+(n-k+1)(k-1)=6n-20 \geq 4n-8.$ Equality
holds if and only if $n=6$, $G$ is a graph composed of a cycle on 4
vertices and a tree on 3 vertices sharing a single vertex.

If $k=4$, then $4n-4k-4+(n-k+1)(k-1)=7n-29 \geq 4n-8.$ Equality
holds if and only if $n=7$, $G$ is a graph composed of a cycle on 4
vertices and a tree on 4 vertices sharing a single vertex.
\end{pf}
\begin{qed}
\end{qed}

Since $G$ is a bipartite graph, $n_0(e)=0$, and thus
$Sz^*(G)=Sz(G)$. So we have the following corollary.

\begin{cor}\label{thm2} [Conjecture \ref{conj}]
Let $G$ be a bipartite connected graph with $n \geq 4$ vertices and
$m \geq n$ edges. Then
$$
Sz^*(G)-W(G) \geq 4n-8.
$$
Moreover the bound is best possible as shown by the graph composed
of a cycle on $4$ vertices $C_4$ and a tree $T$ on $n-3$ vertices
sharing a single vertex.
\end{cor}

\end{document}